\documentclass{article}

\newtheorem{remark}{\bf Remark}[section]
\newtheorem{definition}{\bf Definition}[section]
\usepackage{tensor}
\usepackage{graphicx}
\usepackage[utf8]{inputenc}
\usepackage{amsmath, amssymb}
\usepackage[]{graphicx, url}
\usepackage[]{indentfirst}
\usepackage{siunitx} 
\usepackage{authblk}

\numberwithin{equation}{section}

\title{Compressible Euler Equations on a Sphere and Elliptic-Hyperbolic Property}

\author[1]{Ian Holloway }
\author[2]{Sivaguru S. Sritharan}

\affil[1]{Department of Mathematics, Wright State University, Dayton, OH 45435\\ email: iancholloway@gmail.com}
\affil[2]{M.S. Ramaiah University of Applied Sciences, Bengaluru, India  \\ email: provostsritharan@gmail.com}

\begin{document}
\maketitle
%%%% Article title to be placed here

%%%% Abstract text to be placed here %%%%%%%%%%%%
\begin{abstract}
In this work we systematically derive the governing equations of supersonic conical flow by projecting the 3D Euler equations onto the unit sphere. These equations result from taking the assumption of conical invariance on the 3D flow field. Under this assumption, the compressible Euler equations reduce to a system defined on the surface of the unit sphere. This compressible flow problem has been successfully used to study steady supersonic flow past cones of arbitrary cross section by reducing the number of spatial dimensions from 3 down to 2 while still capturing many of the relevant 3D effects. In this paper the powerful machinery of tensor calculus is utilized to avoid reference to any particular coordinate system. With the flexibility to use any coordinate system on the surface of a sphere, the equations can be more readily solved numerically when a structured mesh is used by defining the mesh lines to be the coordinate lines. The type of the system of partial differential equations would be hyperbolic or elliptic based on whether the crossflow Mach number is supersonic or subsonic.
\end{abstract}
%%%%%%%%%%%%%%%%%%%%%%%%%%%

%%%%%%%%%% Insert the texts which can accomdate on firstpage in the tag "fmtext" %%%%%

\section{Introduction}

Supersonic flow past a cone has been widely studied as a type of 3D flow past a body that is more easily analyzed and provides valuable insight into how flow past a true aircraft will behave. A common assumption to make for for studying fluid flows past bodies is that the flow is uniform in one coordinate direction. In the case of flow past airfoils, the assumption is that the flow is uniform along the length of the wing. Such an assumption reduces the problem from a three dimensional problem to a two dimensional problem defined on a single cross section of the wing. The math involved is simplified and numerical solutions can be achieved more quickly. The conical assumption provides the same benefits of dimension reduction, but provides a different perspective than planar reductions such as the wing cross section. Whereas the wing cross section reduction is a side view, the conical assumption is like a rear view. Studying a 2D cross section of a wing provides insight into flow behavior at the leading and trailing edge, and along the top and bottom surfaces, but it does not include information about flow at the wing tip, or wing root, nor does it provide information about crossflow on the top and bottom surfaces of the wing. In contrast, the conical assumption does not model the leading or trailing edges, but provides insight into flow along the width of a wing, including tip, root, and crossflow. It thus serves to fill in gaps left by other common assumptions.

Maslen  \cite{Maslen} gives a brief history of work done on the subject prior to his paper, consisting largely of linear or thin body approximations. Work on the nonlinear problem was most famously done by Taylor and Maccoll \cite{TaylorMaccoll} who considered a circular cone at zero angle of attack. In \cite{sriThesis}, Sritharan used the machinery of tensor calculus to project the mass equation onto the unit sphere for the case of potential flow. In related works, \cite{FAmethod,Guan}, the more sophisticated finite volume methods available at the time were used to compute numerical solutions for cones of arbitrary cross section.

%history of conical flow

Now, with even more sophisticated numerical methods for fluid flow equations and a renewed interest in hypersonics, we add to this body of knowledge by extending Sritharan's work and deriving the conical form of the full system of Euler equations. The resulting system is not restricted to potential or isentropic flow or any other type of further approximation and is stated for a general coordinate system. This work is done in tandem with \cite{ConMHD} and \cite{ConSolver}. The work  in \cite{ConMHD} is the same analysis done in this paper, but applied to the Ideal Magnetohydrodynamics equations. These equations describe an Euler flow which is electrically conducting such as those that occur in high altitude hypersonic flights. Therefore, there is the same motivation for the conical assumption to be applied. The resulting systems are similar enough that it was convenient to develop a numerical scheme which solves both systems. This scheme is described in \cite{ConSolver}.

%maybe leave these out?
For reference, the standard Euler equations \cite{anderson,liepmann} are presented in Equation \eqref{Euler}.

\begin{subequations}\label{Euler} %might need a reference for these equations - like Serrin
\begin{gather}
      \rho_t  + \nabla\cdot(\rho \pmb{V})=0
\\
      (\rho \pmb{V})_t  + \nabla\cdot(\rho \pmb{V}\otimes\pmb{V}+PI)= \pmb{0}
\\
      (\rho E)_t  + \nabla\cdot(\rho E+P)\pmb{V} = 0 
\end{gather}
\end{subequations}

The dependent variables are $\rho$, $\pmb{V}$, and $e$, and $P=P(\rho,e)$ is provided by a gas law to close the system. To use the machinery of tensor calculus, it is convenient to consider Equation \eqref{Euler} in coordinate free form given in Equation \eqref{transient} for the contravariant components \cite{Eisemann}.

\begin{subequations}\label{transient} 
\begin{gather}
      \left(\sqrt{G}\rho\right)_t  + \left(\rho\sqrt{G}V^j\right)_{|j}=0
\\
      \left(\rho\sqrt{G}V^i\right)_t  + \left(\sqrt{G}\left[\rho V^iV^j + G^{ij}P\right] \right)_{|j}= \pmb{0}
\\
      \left(\rho\sqrt{G}E\right)_t  + \left( \sqrt{G}\left[\rho E+P\right] V^j \right)_{|j}= 0
\end{gather}
\end{subequations}

The notation $(\cdot)_{|j}$ refers to the covariant derivative. The steady problem is the object of consideration (and furthermore time dependency is incompatible with the conical assumption), so the time derivative terms will all be set to zero, leaving Equation \eqref{steady}:

\begin{subequations}\label{steady}
\begin{gather}
      \left(\rho\sqrt{G}V^j\right)_{|j}=0 
\\
      \left(\sqrt{G}\left[\rho V^iV^j + G^{ij}P\right] \right)_{|j}= \pmb{0}
\\
     \left( \sqrt{G}\left[\rho E+P\right] V^j \right)_{|j}= 0\
\end{gather}
\end{subequations}

After applying the conical assumption and projecting the system onto the surface of a sphere, the result is:

\begin{subequations}\label{TheEq}
\begin{gather}
\frac{\partial}{\partial \xi^\beta}\left(\rho\sqrt{g}v^\beta\right) + 2\rho \sqrt{g}V^3 = 0
\label{mass} \\
\frac{\partial }{\partial\xi^\beta}(\sqrt{g}\left[\rho v^\alpha v^\beta + g^{\alpha\beta}P\right]) 
+ \overset{(g)}{\Gamma}\indices{_\gamma^\alpha_\nu}\sqrt{g}\left[\rho v^\gamma v^\nu+g^{\gamma\nu}P\right] 
+ 3\rho\sqrt{g}v^\alpha V^3 = 0
\label{mom} \\
v^\alpha\frac{\partial V^3}{\partial \xi^\alpha} - q^2_c = 0
\label{mom_r} \\
\frac{\partial}{\partial \xi^\beta}\left( \sqrt{g}\left[\rho E+P\right] v^\beta \right) + 2\sqrt{g}\left[\rho E+P\right] V^3 =0 \label{energy}
\end{gather}
\end{subequations}

In all of these equations, Einstein summation is used where summations are carried out over repeated indices (in this article the convention is adopted that Latin indices such as $i,j$ take on values from 1 to 3 and Greek indices such as $\alpha,\beta$ take on values from 1 to 2). Variables are as follows: $\rho$ is the density, $v^\beta$ are velocity components on the surface of a sphere scaled to have no ``$r$'' dependency, $V^3$ is the radial component of velocity, $E=e+\frac{1}{2}|\pmb{V}|^2$ is the total specific energy (thermal plus kinetic), where $e$ is the specific thermal energy, $P=P(\rho, e)$ is the pressure, $g_{\beta\alpha}$ is the metric tensor characterizing angle and distance on the surface of the sphere (with ``$r$'' dependency removed) which has inverse $g^{\omega\alpha}$, $g$ is the determinant of the metric tensor, $\xi^\beta$ are the coordinates on the surface of the sphere, $q_c=\sqrt{g_{\beta\alpha}v^\beta v^\alpha}$ is the magnitude of the crossflow velocity, and $\overset{(g)}{\Gamma}\indices{_\gamma^\alpha_\nu}$ is the Christoffel symbol defined in terms of the metric tensor on the surface of the sphere. These equations make up a new first order system of partial differential equations which has not yet been derived or analyzed by the academic community. 

Defining the speed of sound, $c$, as:

\begin{equation}
    c = \frac{\sqrt{PP_e + \rho^2 P_\rho} }{\rho}
\end{equation}

The type of the system was determined to be:

\begin{equation}
\text{type}=\begin{cases}
	\text{hyperbolic} & q_c > c \\ 
    \text{elliptic} & q_c < c \\ 
\end{cases}
\end{equation}

with characteristic speeds:

\begin{equation}\label{FullEvalsTop}
\frac{v^2}{v^1}, \frac{v^2}{v^1}, \frac{v^2}{v^1},\\
\frac{v^1v^2 - c^2g^{12}\mp \frac{c}{\sqrt{g}}\sqrt{q_c^2-c^2}}{(v^1)^2-g^{11}c^2}
\end{equation}

treating $\xi^1$ as the time-like direction.

In section \ref{ProblemSetting} the setting of the problem is described qualitatively along with key features distinctive of the solution. Section \ref{prelim} describes the geometric machinery necessary to derive system \eqref{TheEq}. Section \ref{conical} describes the conical assumption imposed upon the unknowns. Following that, the projected equations are derived in section \ref{derivation}. The derived equations are compared to the case of using spherical coordinates in section \ref{sphcoords} to demonstrate the consistency. Sections \ref{hyperbolic} and \ref{valsvects} discuss the type of system \eqref{TheEq} based on its eigenvalues, and in section \ref{potential} these results are compared to the case of potential flow which has been studied previously.

\section{Problem Setting}\label{ProblemSetting}

\begin{figure}[ht]
    \centering
    \includegraphics{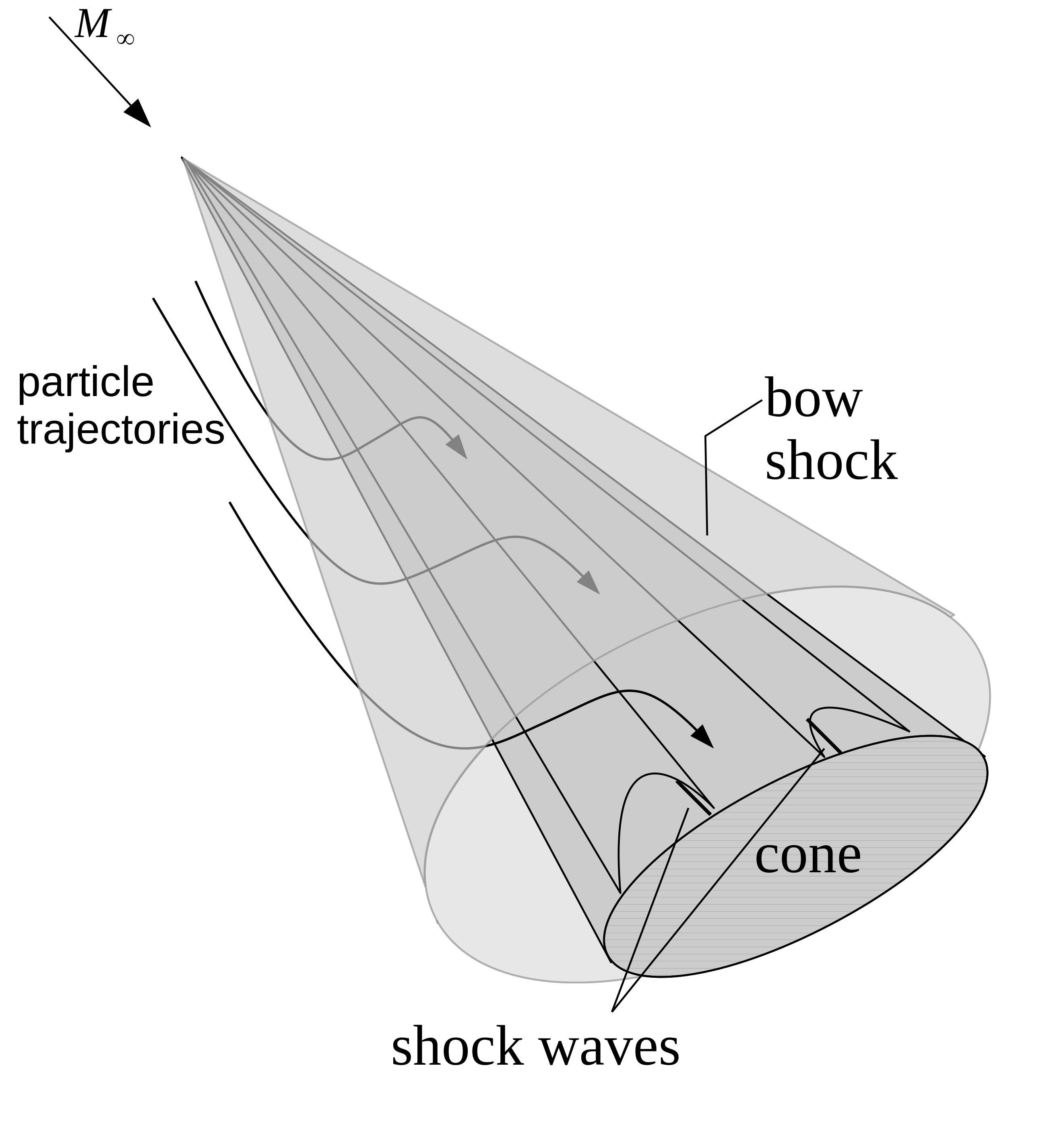}
    \caption{Supersonic infinite cones with elliptic cross section. Shock wave formation and particle trajectories are shown.}
    \label{crossShocks}
\end{figure}

The cone of arbirtrary cross section is considered to be infinite and at an angle of attack relative to the free stream. Features of this flow include an attached bow shock wave, crossflow streamlines which wrap around the body and converge above, and two or more body shocks which are caused by the crossflow briefly going supersonic \cite{sriThesis,ShockFreeCrossFlow,NASA_con,RemarksConFlow}. These various features are depicted in Figure \ref{crossShocks}.

\begin{figure}[ht]
    \centering
    \includegraphics{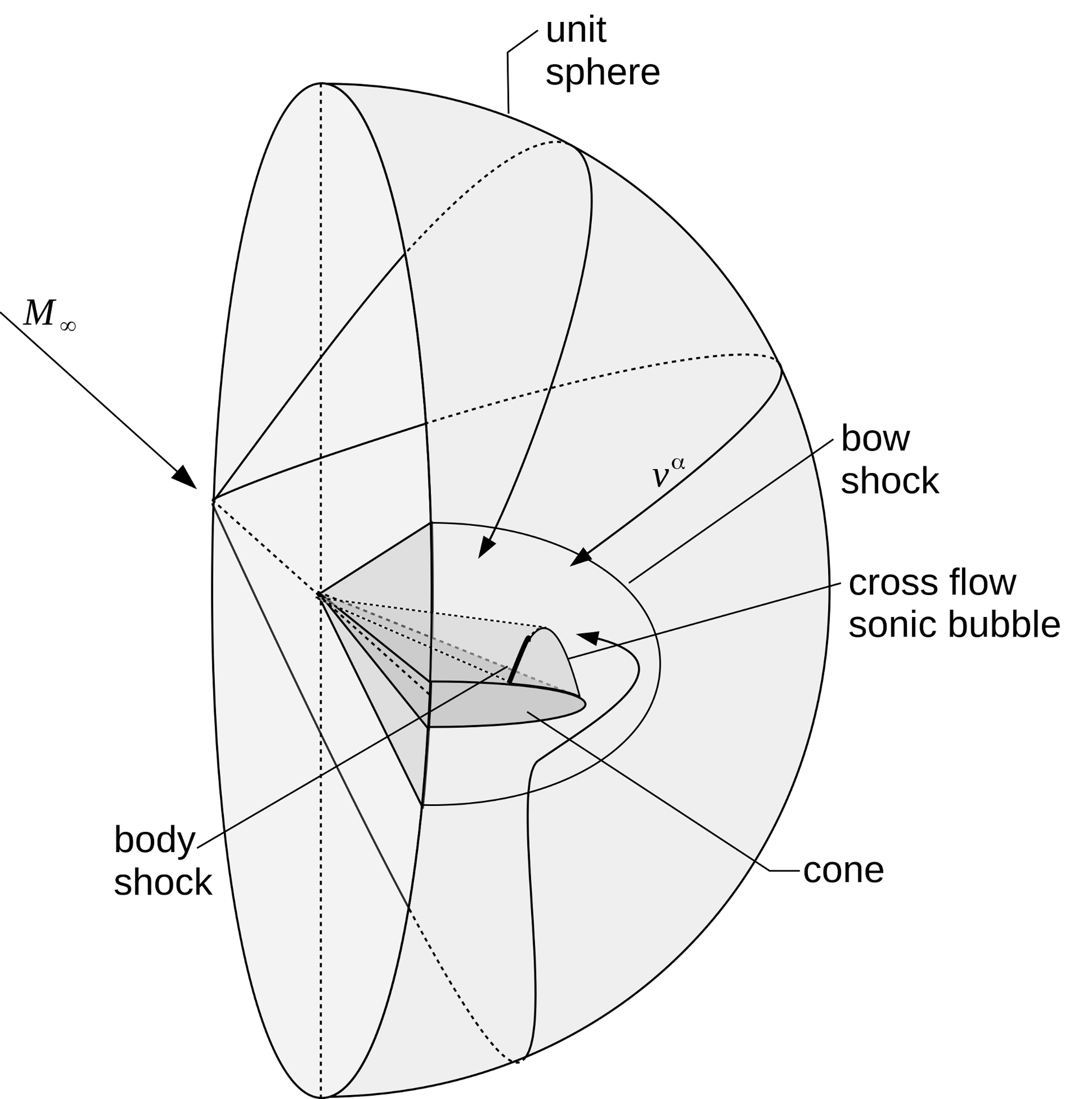}
    \caption{Problem setting sliced by a sphere with the velocity projected onto the surface giving the crossflow streamlines.}
    \label{sphere}
\end{figure}

Such a flow is said to be conical if there exists a point in the domain such that along any line that goes through this point, the flow properties (density, velocity, energy, etc) do not change \cite{sriThesis,RemarksConFlow}. Effectively, this means that if the origin is set to be the tip of the cone, then the solution has no ``$r$'' dependency, where $r$ is the distance from the origin. This type of flow can best be studied by taking a spherical slice out of the domain centered on the origin and projecting the velocity onto that sphere as shown in Figure \ref{sphere}. A solution obtained on this spherical shell of a given radius will thus be valid on a shell of any other radius so that the flow in the whole of the 3D domain is accounted for.

\begin{figure}[ht]
    \centering
    \includegraphics{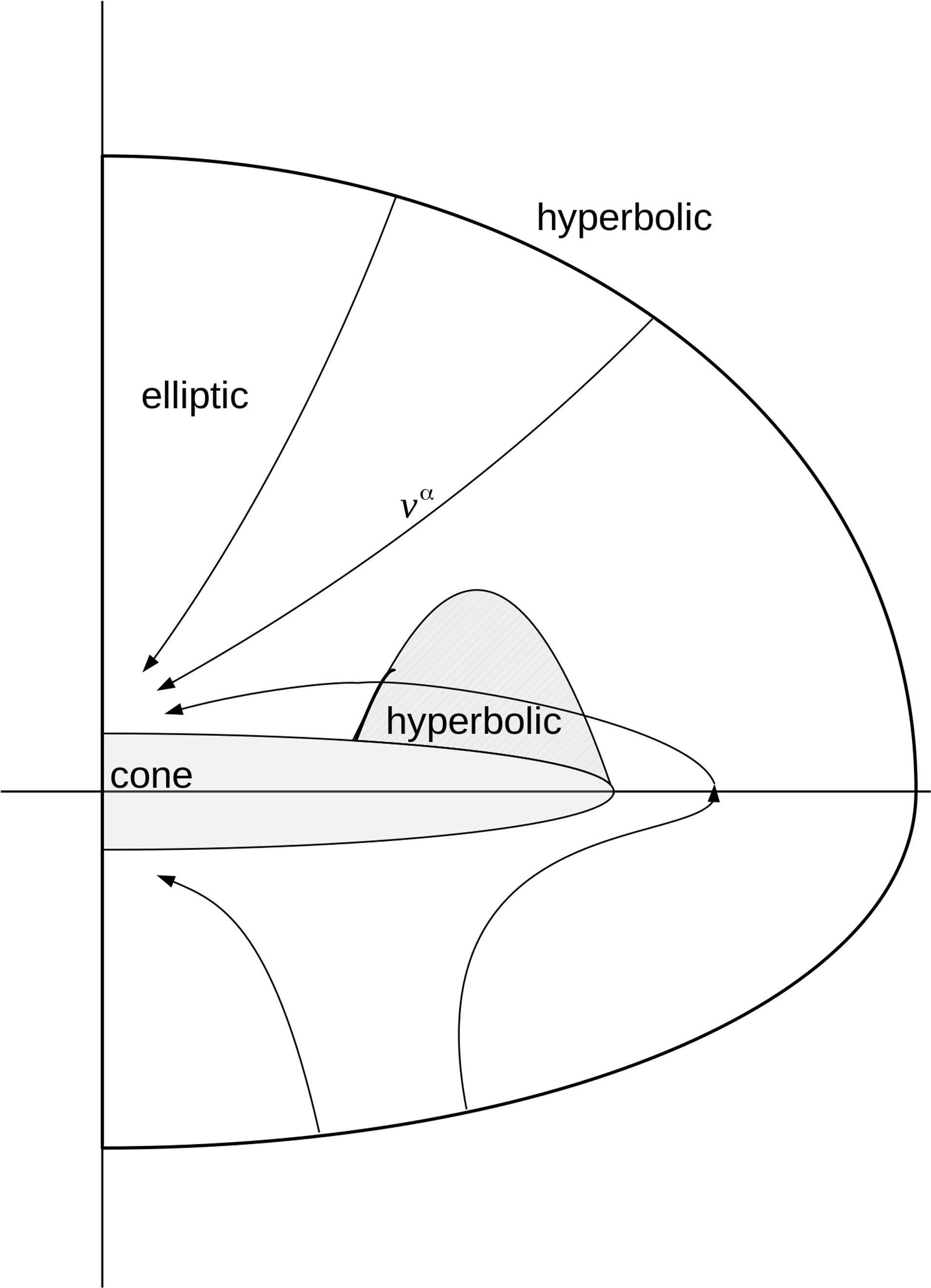}
    \caption{Diagram of types of the governing PDE system. The darkly colored line around the outside is the bow shock. The other darkly colored line on the inner boundary of the hyperbolic bubble is the body shock. }
    \label{rearView}
\end{figure}

Another interesting feature of this flow is that the governing system of partial differential equations changes type multiple times within the domain. As stated previously, the system is hyperbolic when the crossflow is supersonic, and elliptic when it is subsonic. The crossflow free stream is established to be supersonic, but will become subsonic after passing through the bow shock. Surrounded by the post-shock elliptic region is the region where the flow wrapping around the body briefly goes supersonic \cite{sriThesis,Guan}. These regions are diagrammed in Figure \ref{rearView} which is a rear view of the cone.

The changing back and forth of the type throughout the domain as well as regions of different types sharing boundaries must be accounted for in the theory and numerical solving of the governing equations.

\section{Geometric Preliminaries}\label{prelim}

Consider a 3D Euclidean space characterized by metric tensor $G_{ij}$ and coordinates $x^i$. Embedded in this 3D space is a 2D spherical subspace characterized by the metric tensor $\tilde{g}_{\alpha\beta}$ and coordinates $\xi^\alpha$. For such a subspace there is the relationship:

\begin{equation}
    \tilde{g}_{\alpha\beta} = G_{ij}B^i_\alpha B^j_\beta
\end{equation}

where the projection factors are given by:

\begin{equation}
    B^i_\alpha = \frac{\partial x^i}{\partial \xi^\alpha}
\end{equation}

and

\begin{equation}
    B^\alpha_i = \tilde{g}^{\alpha\beta}G_{ij}B^j_\beta
\end{equation}

A tensor in the embedding space can be projected onto the sphere using the projection factors, such as:

\begin{equation}
    \tilde{w}^\alpha = B_i^\alpha W^i, \phantom{m} \tilde{w}^{\alpha\beta} = B_i^\alpha B_j^\beta W^{ij}, \phantom{m} \tilde{w}^{\alpha\beta}_\nu = B_i^\alpha B_j^\beta B^k_\nu W^{ij}_k, \phantom{m} \text{etc.}
\end{equation}

It is convenient to treat the three dimensional embedding space as having the two subspace coordinates and a radial coordinate as its three coordinates. That is $\pmb{x}=(\xi^1,\xi^2,r)$. The $r$ coordinate is orthogonal to the other coordinates so the metric tensor of the embedding space in matrix form would be:

\begin{equation}
    G_{ij} = \left[ \begin{smallmatrix} \cdot&\cdot&0\\ \cdot&\cdot&0 \\ 0&0&1 \end{smallmatrix} \right], \phantom{m} 1\leq i,j \leq 3
\end{equation}

and that of the embedded subspace:

\begin{equation}
    \tilde{g}_{\alpha\beta}= G_{\alpha\beta}, \phantom{m} 1\leq\alpha,\beta\leq 2
\end{equation}

\begin{figure}[ht]
    \centering
    \includegraphics{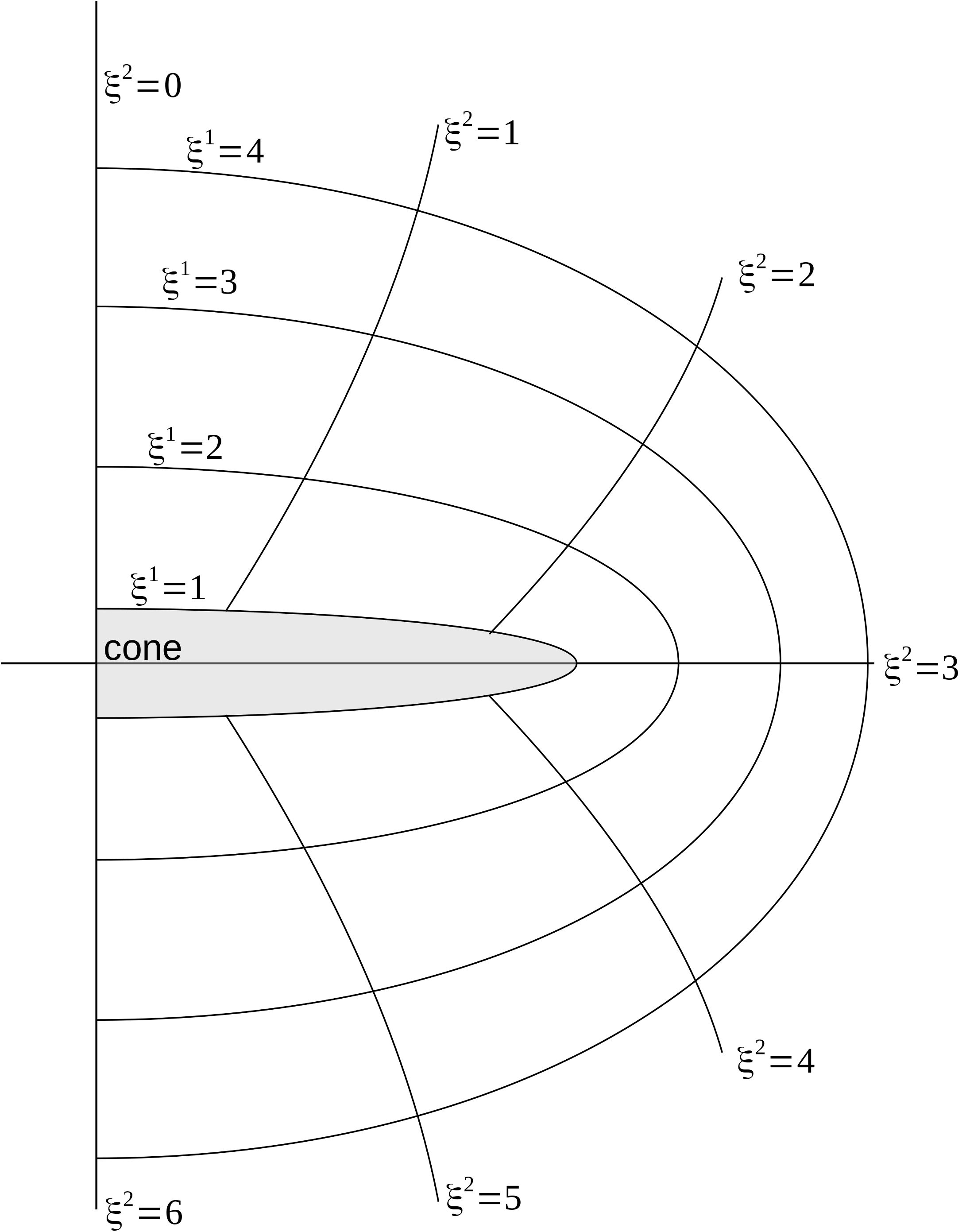}
    \caption{Example of coordinate lines which conform to the shape of the body and are not necessarily orthogonal. }
    \label{sampCoords}
\end{figure}

\begin{remark}
Note that though traditional spherical coordinates $\theta$ and $\phi$ on the surface of the sphere would be a valid choice of coordinates, one is not restricted to them. For this topic, one can consider any two surface coordinates and a radial one. This allows for the possibility of the coordinate lines being aligned with the surface of the cone (as shown in Figure \ref{sampCoords}) even if it has an irregular cross section. In the case of a numerical solution using a structured mesh, the coordinate lines can be defined to follow the mesh lines and simplify some calculations. In particular, it is not necessary to compute dot products with the normal of the computational cell boundary when computing the flux through that interface. 
\end{remark}

Any vector $\tilde{w}^\alpha$ defined at a point in the subspace will have a length defined by:

\begin{equation}
    |\tilde{\pmb{w}}|^2=\tilde{g}_{\alpha\beta}\tilde{w}^\alpha\tilde{w}^\beta
\end{equation}

Since this subspace is defined to be the surface of a sphere, distances will scale proportional to the radius of the subspace, $r$, giving:

\begin{equation}
    |\tilde{\pmb{w}}|^2=\tilde{g}_{\alpha\beta}\tilde{w}^\alpha\tilde{w}^\beta = r^2g_{\alpha\beta}\tilde{w}^\alpha\tilde{w}^\beta
\end{equation}

where the $r$ dependency has been separated out of the metric tensor. This implies that $\tilde{g}_{\alpha\beta}=r^2g_{\alpha\beta}$ (also $\tilde{g}^{\alpha\beta}=\frac{1}{r^2}g^{\alpha\beta}$) and that $g_{\alpha\beta}$ is a function of $\xi^1$ and $\xi^2$ only. This leads us to define a new representation of the vector where $w^\alpha=r\tilde{w}^\alpha$ and also $w_\alpha=\frac{1}{r}\tilde{w}_\alpha$. Using this definition:

\begin{equation}\label{rescaledLength}
    |\tilde{\pmb{w}}|^2=\tilde{g}_{\alpha\beta}\tilde{w}^\alpha\tilde{w}^\beta = r^2g_{\alpha\beta}\tilde{w}^\alpha\tilde{w}^\beta = g_{\alpha\beta}w^\alpha w^\beta
\end{equation}

In particular, equation \eqref{rescaledLength} says the magnitude of the surface components of a vector does not change as you scale in r. This representation, with the $r$ dependency shifted from the metric onto the vector components, can be used for any vector.

\section{Conical assumption}\label{conical}

This article describes the Euler equations subject to the conical assumption on all the dependent variables (density, velocity, and energy). 

\begin{definition}
A quantity is said to be conical if the covariant derivative in the $r$ direction is identically zero.
\end{definition}

For scalar quantities such as $\rho$ and $E$, this means that the partial derivative with respect to $r$ is zero. For higher order tensorial quantities it is not so simple. Because the basis for the vectors is not uniform, it is possible for the components of a vector to change, but for the vector to remain the same, and conversely for the vector to change, but the components to remain the same. Therefore the covariant derivative must be used, which accounts for the changing of the underlying coordinate basis. 

Consider a vector, $\pmb{W}$, in the 3D embedding space. It has 3 components; two corresponding to the spherical subspace and one radial component, that is $\pmb{W}=\left[ \begin{smallmatrix} \tilde{w}^1 & \tilde{w}^2 & W^3 \end{smallmatrix} \right]$. If $\pmb{W}$ is conical, then all components of the covariant derivative in the $x^3$ (or $r$) direction are identically zero. Mathematically, that is:

\begin{equation}
W^i_{|3}=0, \phantom{m}\forall i
\end{equation}

Inserting the full expression for the covariant derivative gives:

\begin{equation*}
    W^i_{|3} = \frac{\partial W^i}{\partial x^3} + \Gamma\indices{_j^i_3}W^j
\end{equation*}

\begin{equation*}
   = \frac{\partial W^i}{\partial x^3} + \frac{G^{ik}}{2}\left[ 
\frac{\partial G_{kj}}{\partial x^3} + \frac{\partial G_{k3}}{\partial x^j} - \frac{\partial G_{j3}}{\partial x^k} \right]W^j
\end{equation*}

Because of the form of the metric, the last two terms in the Christoffel symbol are identically zero, giving:

\begin{equation}
   = \frac{\partial W^i}{\partial x^3} + \frac{G^{ik}}{2}\left[ 
\frac{\partial G_{kj}}{\partial x^3} \right]W^j
\end{equation}

Examining this expression, when $i=3$ it becomes:

\begin{equation}
    \frac{\partial W^3}{\partial x^3} = \frac{\partial W^3}{\partial r} =0
\end{equation}

And otherwise:

\begin{equation}
   = \frac{\partial \tilde{w}^\alpha}{\partial x^3} + \frac{\tilde{g}^{\alpha\nu}}{2}\left[ 
\frac{\partial \tilde{g}_{\nu\beta}}{\partial x^3} \right]\tilde{w}^\beta
\end{equation}

Plugging in the components with the shifted $r$ dependency gives:

\begin{equation*}
   = \frac{\partial }{\partial r}\left(\frac{1}{r}w^\alpha\right) + \frac{g^{\alpha\nu}}{2r^2}\left[ 
\frac{\partial }{\partial r}\left(r^2g_{\nu\beta}\right) \right]\frac{1}{r}w^\beta
\end{equation*}

\begin{equation*}
   = \frac{1}{r}\frac{\partial }{\partial r}\left(w^\alpha\right) - \frac{1}{r^2}w^\alpha + \frac{g^{\alpha\nu}}{2r^2}\left[ 
2rg_{\nu\beta} \right]\frac{1}{r}w^\beta
\end{equation*}

\begin{equation*}
   = \frac{1}{r}\frac{\partial w^\alpha}{\partial r} - \frac{1}{r^2}w^\alpha + \frac{1}{r^2}\delta^{\alpha}_\beta w^\beta
\end{equation*}

\begin{equation*}
   = \frac{1}{r}\frac{\partial w^\alpha}{\partial r} - \frac{1}{r^2}w^\alpha + \frac{1}{r^2}w^\alpha
\end{equation*}

\begin{equation}
   = \frac{1}{r}\frac{\partial w^\alpha}{\partial r}=0
\end{equation}

Thus the conical assumption implies: 

\begin{equation}
   \frac{\partial W^3}{\partial r} = \frac{\partial w^\alpha}{\partial r} = 0
\end{equation}

This expression tells us that it is the rescaled components of the vector and not the original components which are independent of $r$. This is an important concept to keep in mind as the conical Euler equations are derived.

\section{Derivation of conical equations}\label{derivation}

For the projection of the equations, the following relations are necessary which involve the various elements of the different spaces:

\begin{equation}\label{metric1}
g_{\alpha\beta} = \frac{1}{r^2}\tilde{g}_{\alpha\beta} = \frac{1}{r^2}G_{\alpha\beta},\phantom{m}
	G_{ij} = \left[ \begin{smallmatrix} \cdot&\cdot&0\\ \cdot&\cdot&0 \\ 0&0&1 \end{smallmatrix} \right]
\end{equation}

\begin{equation}\label{metric2}
g^{\alpha\beta} = r^2\tilde{g}^{\alpha\beta} 
\end{equation}

\begin{equation}
g = \frac{1}{r^4}\tilde{g} = \frac{1}{r^4}G
\end{equation}

\begin{equation}
\sqrt{g} = \frac{1}{r^2}\sqrt{\tilde{g}} = \frac{1}{r^2}\sqrt{G}
\end{equation}

\begin{equation}
\tilde{g}_{\alpha\beta} = G_{ij}B^i_\alpha B^j_\beta
\end{equation}

\begin{equation}
\tilde{g}^{\alpha\beta}B^j_\beta = G^{ij}B_i^\alpha 
\end{equation}

\begin{equation}\label{projRank1}
\tilde{v}^\alpha = B^\alpha_iV^i,\phantom{m} B^i_\alpha = \frac{\partial x^i}{\partial\xi^\alpha}
\end{equation}

\begin{equation}
    B^\beta_jB^i_\beta = \delta^i_j - N^i_j,\phantom{m} N^i_j \equiv \delta^i_3\delta^3_j
\end{equation}

\begin{equation}
\tilde{v}_\alpha = rv_\alpha
\end{equation}

\begin{equation}
\tilde{v}^\alpha = \frac{1}{r}v^\alpha
\end{equation}

%For a rank 1 contravariant relative tensor of weight 1, we have:
%
%\begin{equation}
%w^{\alpha} = \rho\sqrt{g}v^\alpha = \rho\frac{1}{r^2}\sqrt{\tilde{g}}r\tilde{v}^\alpha = \frac{1}{r}\tilde{w}^{\alpha\beta} \\
%= \frac{1}{r}B^\alpha_i\rho\sqrt{G}V^i = \frac{1}{r}B^\alpha_iW^{i}
%\end{equation}

%For a rank 2 contravariant relative tensor of weight 1, we have:

%\begin{equation}
%w^{\alpha\beta} = \rho\sqrt{g}v^\alpha v^\beta = \rho\sqrt{\tilde{g}}\tilde{v}^\alpha \tilde{v}^\beta = \tilde{w}^{\alpha\beta} \\
%= B^\alpha_iB^\beta_j\rho\sqrt{G}V^iV^j = B^\alpha_iB^\beta_jW^{ij}
%\end{equation}

The projection begins by first considering the continuity and energy equations. For these equations, the LHS's are the contracted covariant derivative (divergence) of rank 1 relative tensors of weight 1 as defined by \cite{Lovelock}. This means that it carries with it the square root of the determinant of the metric tensor raised to the first power, and that when it transforms from one coordinate system to the other, it changes metric determinants as well. The contracted covariant derivative of these tensors is given by the following:

Let $W\indices{^{j}}$ be a rank 1 contravariant relative tensor of weight 1, such as $\rho\sqrt{G}V^j$ or $\sqrt{G}\left[\rho E+P\right] V^j$. Then the contracted covariant derivative is given by:

\begin{equation}
W\indices{^{j}_{|j}} = \frac{\partial W^{j}}{\partial x^j} 
\end{equation}

And likewise for a surface tensor of the same type and weight, and its rescaled components:

\begin{equation}
\tilde{w}\indices{^{\beta}_{||\beta}} = \frac{\partial \tilde{w}^{\beta}}{\partial \xi^\beta} 
\end{equation}

\begin{equation}
w\indices{^{\beta}_{||\beta}} = \frac{\partial w^{\beta}}{\partial \xi^\beta} 
\end{equation}

In this case $\tilde{w}^\beta=\rho\sqrt{\tilde{g}}\tilde{v}^\beta$ or $\sqrt{\tilde{g}}[\rho E + P]\tilde{v}^\beta$ and $w^\beta=\rho\sqrt{g}v^\beta$ or $\sqrt{g}[\rho E + P]v^\beta$.

%BEFORE THIS POINT MUST HAVE EXPLAINED RESCALING RELATIVE TENSORS
Using previous relations, the surface divergence for the rescaled components can be rewritten as:

\begin{equation*}
w\indices{^{\beta}_{||\beta}} = \frac{\partial }{\partial \xi^\beta}\left(\frac{1}{r}\tilde{w}^\beta\right) = \frac{1}{r}\frac{\partial }{\partial \xi^\beta}\left(B^\beta_iW^i\right) 
\end{equation*}

\begin{equation*}
= \frac{1}{r}B^\beta_i\frac{\partial W^i}{\partial \xi^\beta} + \frac{1}{r}W^i\frac{\partial }{\partial \xi^\beta}B^\beta_i
\end{equation*}

\begin{equation*}
= \frac{1}{r}B^\beta_iB^j_\beta\frac{\partial W^i}{\partial x^j} + 0
\end{equation*}

\begin{equation*}
= \frac{1}{r}\left(\delta^i_j - N^i_j\right)\frac{\partial W^i}{\partial x^j}
\end{equation*}

\begin{equation*}
= \frac{1}{r}\frac{\partial W^j}{\partial x^j} - \frac{1}{r}\frac{\partial W^3}{\partial x^3}
\end{equation*}

Thus:

\begin{equation}\label{projDivRank1}
w\indices{^{\beta}_{||\beta}}= \frac{1}{r}W\indices{^j_{|j}} - \frac{1}{r}\frac{\partial W^3}{\partial x^3}
\end{equation}

The mass and energy flux vectors can now be plugged into this expression. For both cases the first term on the RHS of Equation \eqref{projDivRank1} is equal to zero so the result is:

\begin{equation*}
\left(\rho\sqrt{g}v^\beta\right)_{||\beta}=-\frac{1}{r}\frac{\partial}{\partial r}\left(\rho\sqrt{G}V^3\right) 
\end{equation*}

\begin{equation*}
\left( \sqrt{g}\left[\rho E+P\right] v^\beta \right)_{||\beta}=-\frac{1}{r}\frac{\partial}{\partial r}\left( \sqrt{G}\left[\rho E+P\right] V^3 \right) 
\end{equation*}

Putting the RHS's into the rescaled form:

\begin{equation*}
\left(\rho\sqrt{g}v^\beta\right)_{||\beta}=-\frac{1}{r}\frac{\partial}{\partial r}\left(\rho r^2\sqrt{g}V^3\right) 
\end{equation*}

\begin{equation*}
\left( \sqrt{g}\left[\rho E+P\right] v^\beta \right)_{||\beta}=-\frac{1}{r}\frac{\partial}{\partial r}\left( r^2\sqrt{g}\left[\rho E+P\right] V^3 \right) 
\end{equation*}

And so:

\begin{equation*}
\left(\rho\sqrt{g}v^\beta\right)_{||\beta}=-2\rho \sqrt{g}V^3 
\end{equation*}

\begin{equation*}
\left( \sqrt{g}\left[\rho E+P\right] v^\beta \right)_{||\beta}=-2\sqrt{g}\left[\rho E+P\right] V^3
\end{equation*}

Finally:

\begin{equation}
\frac{\partial}{\partial \xi^\beta}\left(\rho\sqrt{g}v^\beta\right) + 2\rho \sqrt{g}V^3 = 0
\end{equation}

\begin{equation}
\frac{\partial}{\partial \xi^\beta}\left( \sqrt{g}\left[\rho E+P\right] v^\beta \right) + 2\sqrt{g}\left[\rho E+P\right] V^3 = 0
\end{equation}

\subsection{Momentum equation}

The projection continues by considering the momentum equation. The LHS is the contracted covariant derivative (divergence) of a rank 2 relative tensor of weight 1 as defined by \cite{Lovelock}. The contracted covariant derivative of which is given by the following:

Let $W\indices{^{ij}}$ be a rank 2 contravariant relative tensor of weight 1. Then the contracted covariant derivative is given by:

\begin{equation}
    W\indices{^{ij}_{|j}} = \frac{\partial W^{ij}}{\partial x^j} + \overset{(G)}{\Gamma}\indices{_h^i_k}W^{hk}
\end{equation}

And likewise for a surface tensor of the same type and weight and its rescaled components.

\begin{equation}
    \tilde{w}\indices{^{\alpha\beta}_{||\beta}} = \frac{\partial \tilde{w}^{\alpha\beta}}{\partial \xi^\beta} + \overset{(\tilde{g})}{\Gamma}\indices{_\gamma^\alpha_\nu}\tilde{w}^{\gamma\nu}
\end{equation}

\begin{equation}\label{surfDivRank2}
    w\indices{^{\alpha\beta}_{||\beta}} = \frac{\partial w^{\alpha\beta}}{\partial \xi^\beta} + \overset{(g)}{\Gamma}\indices{_\gamma^\alpha_\nu}w^{\gamma\nu}
\end{equation}

Where the Christoffel symbols are defined in terms of the respective metric tensors. In analogy to the projection of the continuity and energy equations, previous relations are plugged into the surface divergence expression for the rescaled tensor. For clarity, the Christoffel symbol is considered by itself first.

\subsubsection{Projected Christoffel Symbol}

To begin, an expression for the Christoffel symbol defined by the rescaled metric which has no $r$ dependency is found. It is given by:

\begin{equation*}
\overset{(g)}{\Gamma}\indices{_\gamma^\alpha_\nu} = \frac{ g^{\alpha\beta}}{2}\left[ 
\frac{\partial g_{\beta\gamma}}{\partial\xi^\nu} + \frac{\partial g_{\beta\nu}}{\partial\xi^\gamma} - \frac{\partial g_{\gamma\nu}}{\partial\xi^\beta} \right] = \frac{\tilde{g}^{\alpha\beta}}{2}\left[ 
\frac{\partial\tilde{g}_{\beta\gamma}}{\partial\xi^\nu} + \frac{\partial\tilde{g}_{\beta\nu}}{\partial\xi^\gamma} - \frac{\partial\tilde{g}_{\gamma\nu}}{\partial\xi^\beta} \right] = \overset{(\tilde{g})}{\Gamma}\indices{_\gamma^\alpha_\nu}
\end{equation*}

because the $r^2$'s from Equations \eqref{metric1} and \eqref{metric2} cancel. Proceeding with the projection:

\begin{equation*}
 = \frac{\tilde{g}^{\alpha\beta}}{2}\left[ 
\frac{\partial }{\partial\xi^\nu}(G_{hk}B^h_\beta B^k_\gamma) + \frac{\partial }{\partial\xi^\gamma}(G_{hk}B^h_\beta B^k_\nu) - \frac{\partial }{\partial\xi^\beta}(G_{hk}B^h_\gamma B^k_\nu) \right]
\end{equation*}

Now pulling out the projection factors and changing derivatives:

\begin{equation*}
 = \frac{\tilde{g}^{\alpha\beta}}{2}\left[ 
B^h_\beta B^k_\gamma B^i_\nu\frac{\partial }{\partial x^i}(G_{hk}) + B^h_\beta B^k_\nu B^i_\gamma\frac{\partial }{\partial x^i}(G_{hk}) - B^h_\gamma B^k_\nu B^i_\beta\frac{\partial }{\partial x^i}(G_{hk}) \right]
\end{equation*}

Adjusting indices on the last term:

\begin{equation*}
 = \frac{\tilde{g}^{\alpha\beta}}{2}\left[ 
B^h_\beta B^k_\gamma B^i_\nu\frac{\partial }{\partial x^i}(G_{hk}) + B^h_\beta B^k_\nu B^i_\gamma\frac{\partial }{\partial x^i}(G_{hk}) - B^i_\gamma B^k_\nu B^h_\beta\frac{\partial }{\partial x^h}(G_{ik}) \right]
\end{equation*}

Then pulling out the common projection factor:

\begin{equation*}
 = \frac{\tilde{g}^{\alpha\beta}}{2}B^h_\beta\left[ 
B^k_\gamma B^i_\nu\frac{\partial }{\partial x^i}(G_{hk}) + B^k_\nu B^i_\gamma\frac{\partial }{\partial x^i}(G_{hk}) - B^i_\gamma B^k_\nu\frac{\partial }{\partial x^h}(G_{ik}) \right]
\end{equation*}

\begin{equation*}
 = \frac{G^{hl}}{2}B^\alpha_l\left[ 
B^k_\gamma B^i_\nu\frac{\partial }{\partial x^i}(G_{hk}) + B^k_\nu B^i_\gamma\frac{\partial }{\partial x^i}(G_{hk}) - B^i_\gamma B^k_\nu\frac{\partial }{\partial x^h}(G_{ik}) \right]
\end{equation*}

Adjusting indices on the first term so all projection factors can be pulled out:

\begin{equation*}
 = \frac{G^{hl}}{2}B^\alpha_l\left[ 
B^i_\gamma B^k_\nu\frac{\partial }{\partial x^k}(G_{hi}) + B^k_\nu B^i_\gamma\frac{\partial }{\partial x^i}(G_{hk}) - B^i_\gamma B^k_\nu\frac{\partial }{\partial x^h}(G_{ik}) \right]
\end{equation*}

\begin{equation*}
 = B^\alpha_lB^i_\gamma B^k_\nu\frac{G^{hl}}{2}\left[ 
\frac{\partial }{\partial x^k}(G_{hi}) + \frac{\partial }{\partial x^i}(G_{hk}) - \frac{\partial }{\partial x^h}(G_{ki}) \right]
\end{equation*}

\begin{equation*}
 = B^\alpha_lB^i_\gamma B^k_\nu\overset{(G)}{\Gamma}\indices{_i^l_k}
\end{equation*}

Thus:

\begin{equation}
 \overset{(g)}{\Gamma}\indices{_\gamma^\alpha_\nu} = \overset{(\tilde{g})}{\Gamma}\indices{_\gamma^\alpha_\nu} = B^\alpha_lB^i_\gamma B^k_\nu \overset{(G)}{\Gamma}\indices{_i^l_k}
\end{equation}

Which shows that the Christoffel symbol projects just like a tensor in this case.

\subsubsection{Projected Momentum Equation}

The full momentum equation can now be projected, starting with Equation \eqref{surfDivRank2}:

\begin{equation*}
\frac{\partial w^{\alpha\beta}}{\partial \xi^\beta} + \overset{(g)}{\Gamma}\indices{_\gamma^\alpha_\nu}w^{\gamma\nu} = \frac{\partial \tilde{w}^{\alpha\beta}}{\partial \xi^\beta} + \overset{(\tilde{g})}{\Gamma}\indices{_\gamma^\alpha_\nu}\tilde{w}^{\gamma\nu}
\end{equation*}

\begin{equation*}
= \frac{\partial}{\partial\xi^\beta}(B^\alpha_iB^\beta_jW^{ij}) + B^\gamma_iB^\nu_jB^m_\gamma B^\alpha_lB^k_\nu \overset{(G)}{\Gamma}\indices{_m^l_k} W^{ij}
\end{equation*}

The projection factors in the first term are pulled out and the chain rule is used to change the derivative. Indices l and i in the second term are also switched.

\begin{equation*}
 = B^\alpha_iB^\beta_jB^f_\beta\frac{\partial}{\partial x^f}(W^{ij}) + B^\alpha_iB^\gamma_lB^m_\gamma B^\nu_jB^k_\nu \overset{(G)}{\Gamma}\indices{_m^i_k} W^{lj}
\end{equation*}

\begin{equation*}
 = B^\alpha_i\left[ (\delta^f_j-N^f_j)\frac{\partial}{\partial x^f}(W^{ij}) + (\delta^m_l-N^m_l)(\delta^k_j-N^k_j) \overset{(G)}{\Gamma}\indices{_m^i_k} W^{lj} \right]
\end{equation*}

\begin{equation*}
 = B^\alpha_i\left[ \frac{\partial}{\partial x^j}(W^{ij}) - \frac{\partial}{\partial x^3}(W^{i3}) + (\delta^m_l\delta^k_j-\delta^m_lN^k_j  -N^m_l\delta^k_j+N^m_lN^k_j) \overset{(G)}{\Gamma}\indices{_m^i_k} W^{lj} \right]
\end{equation*}

\begin{equation*}
 = B^\alpha_i\left[ \left(\frac{\partial}{\partial x^j}(W^{ij}) + \overset{(G)}{\Gamma}\indices{_l^i_j}W^{lj}\right) - \frac{\partial}{\partial x^3}(W^{i3}) - \left( \overset{(G)}{\Gamma}\indices{_l^i_3}W^{l3} + \overset{(G)}{\Gamma}\indices{_3^i_j}W^{3j} \right)
+ \overset{(G)}{\Gamma}\indices{_3^i_3}W^{33} \right]
\end{equation*}

The first term in parenthesis is $W\indices{^{ij}_{|j}}$, which for the momentum flux is equal to $0$. The terms in the second set of parenthesis are equal due to symmetry of all the terms. Using this information:

\begin{equation}\label{surfMomentum}
w\indices{^{\alpha\beta}_{||\beta}} = B^\alpha_i\left[ - \frac{\partial}{\partial x^3}(W^{i3}) - 2W^{l3}\overset{(G)}{\Gamma}\indices{_l^i_3} + W^{33}\overset{(G)}{\Gamma}\indices{_3^i_3} \right]
\end{equation}

For clarity, each term is treated individually.

\subsubsection{Term 1}

\begin{equation*}
B^\alpha_i\frac{\partial}{\partial x^3}(W^{i3}) = B^\alpha_i\frac{\partial}{\partial r}(\sqrt{G}\left[\rho V^iV^3 + G^{i3}P\right])
\end{equation*}

\begin{equation*}
= \frac{\partial}{\partial r}(\sqrt{G}\left[\rho\tilde{v}^\alpha V^3 + \tilde{g}^{\alpha\omega}B^3_\omega P\right])
\end{equation*}

\begin{equation*}
= \frac{\partial}{\partial r}(\rho r^2\sqrt{g}\frac{1}{r}v^\alpha V^3)
\end{equation*}

\begin{equation}
= \rho\sqrt{g}v^\alpha V^3
\end{equation}

\subsubsection{Term 2}

\begin{equation*}
B^\alpha_iW^{l3}\overset{(G)}{\Gamma}\indices{_l^i_3} = B^\alpha_iW^{l3}
\frac{G^{ij}}{2}\left[ \frac{\partial }{\partial x^3}(G_{jl}) + \frac{\partial }{\partial x^l}(G_{j3}) - \frac{\partial }{\partial x^j}(G_{l3}) \right]
\end{equation*}

\begin{equation*}
= W^{l3}\frac{\tilde{g}^{\alpha\omega}}{2}B^j_\omega\left[ \frac{\partial }{\partial x^3}(G_{jl}) + \frac{\partial }{\partial x^l}(G_{j3}) - \frac{\partial }{\partial x^j}(G_{l3}) \right]
\end{equation*}

\begin{equation*}
= \sqrt{G}\left[\rho V^lV^3 + G^{l3}P\right]\frac{\tilde{g}^{\alpha\omega}}{2}\left[ B^j_\omega\left(\frac{\partial }{\partial x^3}(G_{jl}) + \frac{\partial }{\partial x^l}(G_{j3})\right) - \frac{\partial }{\partial\xi^\omega}(G_{l3}) \right]
\end{equation*}

Because $B^j_\omega$ acts like $\delta^j_\omega$, and because the $r$ coordinate is orthogonal to the surface coordinates and scales uniformly, this expression becomes:

\begin{equation*}
= \sqrt{G}\left[\rho\tilde{v}^\beta V^3 + G^{\beta 3}P\right]\frac{\tilde{g}^{\alpha\omega}}{2}\left[ \frac{\partial }{\partial x^3}(G_{\omega \beta})  \right]
\end{equation*}

\begin{equation*}
= \sqrt{G}\left[\rho\tilde{v}^\beta V^3\right]\frac{\tilde{g}^{\alpha\omega}}{2}\frac{\partial }{\partial x^3}(\tilde{g}_{\omega \beta})
\end{equation*}

\begin{equation*}
= \rho r^2\sqrt{g}\frac{1}{r}v^\beta V^3\frac{\tilde{g}^{\alpha\omega}}{2}\frac{\partial }{\partial r}(r^2g_{\omega \beta})
\end{equation*}

\begin{equation*}
= \rho r^2\sqrt{g}\frac{1}{r}v^\beta V^3\frac{1}{r^2}\frac{g^{\alpha\omega}}{2}(2rg_{\omega \beta})
\end{equation*}

\begin{equation*}
= \rho \sqrt{g}v^\beta\delta^\alpha_\beta V^3
\end{equation*}

\begin{equation}
= \rho \sqrt{g}v^\alpha V^3
\end{equation}

\subsubsection{Term 3}

\begin{equation}
B^\alpha_iW^{33}\overset{(G)}{\Gamma}\indices{_3^i_3} = B^\alpha_iW^{33}\frac{G^{ij}}{2}\left[ \frac{\partial }{\partial x^3}(G_{j3}) + \frac{\partial }{\partial x^3}(G_{j3}) - \frac{\partial }{\partial x^j}(G_{33}) \right]
\end{equation}

Because the $r$ coordinate is orthogonal to the surface coordinates and scales uniformly, this expression is identically zero.

\subsubsection{Return to Momentum Equation}

Plugging the revised formulas for the 4 terms into Equation \eqref{surfMomentum} gives:

\begin{equation}\label{surfMomRev}
w\indices{^{\alpha\beta}_{||\beta}} = -\rho\sqrt{g}v^\alpha V^3 - 2\rho\sqrt{g}v^\alpha V^3 + 0
\end{equation}

\begin{equation}\label{surfMomFin1}
w\indices{^{\alpha\beta}_{||\beta}} = - 3\rho\sqrt{g}v^\alpha V^3
\end{equation}

Thus:

\begin{equation}\label{surfMomFin2a}
\frac{\partial }{\partial\xi^\beta}(\sqrt{g}\left[\rho v^\alpha v^\beta+g^{\alpha\omega}P\right]) 
+ \overset{(g)}{\Gamma}\indices{_\gamma^\alpha_\nu}\sqrt{g}\left[\rho v^\alpha v^\beta+g^{\alpha\omega}P\right] = 
- 3\rho\sqrt{g}v^\alpha V^3
\end{equation}

or

\begin{equation}\label{surfMomFin2b}
\frac{\partial }{\partial\xi^\beta}(\sqrt{g}\left[\rho v^\alpha v^\beta+g^{\alpha\omega}P\right]) 
+ \overset{(g)}{\Gamma}\indices{_\gamma^\alpha_\nu}\sqrt{g}\left[\rho v^\alpha v^\beta+g^{\alpha\omega}P\right] + 3\rho\sqrt{g}v^\alpha V^3 = 0
\end{equation}

\subsection{Final Equation}

The above equation is really two equations meaning that with the projected mass and energy equations there are 4 total equations but 5 unknowns. The last equation comes from looking at the 3rd component of the momentum equation.

\begin{equation*}
W\indices{^{3j}_{|j}} = 0
\end{equation*}

Expanding out the LHS:

\begin{equation*}
W\indices{^{3j}_{|j}} = \frac{\partial W^{3j}}{\partial x^j} + \overset{(G)}{\Gamma}\indices{_h^3_j}W^{hj}
\end{equation*}

\begin{equation*}
= \frac{\partial}{\partial x^j}(\sqrt{G}\left[\rho V^jV^3+G^{j3}P\right]) + \frac{G^{3l}}{2}\left[ 
\frac{\partial }{\partial x^j}(G_{lh}) + \frac{\partial }{\partial x^h}(G_{lj}) - \frac{\partial }{\partial x^l}(G_{hj}) \right]W^{hj}
\end{equation*}

Separating the surface and radial derivatives:

\begin{multline*}
= \frac{\partial}{\partial \xi^\alpha}(\rho\sqrt{G}\tilde{v}^\alpha V^3) + \frac{\partial}{\partial x^3}(\sqrt{G}\left[\rho(V^3)^2+P\right]) \\ + \frac{G^{3l}}{2}\left[ 
\frac{\partial }{\partial x^j}(G_{lh}) + \frac{\partial }{\partial x^h}(G_{lj}) - \frac{\partial }{\partial x^l}(G_{hj}) \right]W^{hj}
\end{multline*}

Changing metrics and keeping in mind that the radial coordinate is orthogonal to the surface coordinates and has unit scale:

\begin{multline*}
= \frac{\partial}{\partial \xi^\alpha}(\rho r^2\sqrt{g}\frac{1}{r}v^\alpha V^3) + \frac{\partial}{\partial r}(\rho r^2\sqrt{g}(V^3)^2) + \frac{\partial}{\partial r}( r^2\sqrt{g}P) \\ + \frac{1}{2}\left[ 
\frac{\partial }{\partial x^j}(G_{3h}) + \frac{\partial }{\partial x^h}(G_{3j}) - \frac{\partial }{\partial x^3}(G_{hj}) \right]W^{hj}
\end{multline*}

\begin{equation*}
= r\frac{\partial}{\partial \xi^\alpha}(\rho \sqrt{g}V^\alpha V^3) + 2r\rho\sqrt{g}(V^3)^2 + 2r\sqrt{g}P + \frac{1}{2}\left[ 
0 + 0 - \frac{\partial }{\partial r}(\tilde{g}_{\alpha\beta}) \right]\tilde{w}^{\alpha\beta}
\end{equation*}

\begin{equation*}
= r\left[\frac{\partial}{\partial \xi^\alpha}(\rho \sqrt{g}v^\alpha V^3) + 2\rho\sqrt{g}(V^3)^2\right] + 2r\sqrt{g}P - \frac{1}{2}\frac{\partial }{\partial r}(r^2g_{\alpha\beta})\sqrt{\tilde{g}}\left[\rho\tilde{v}^\alpha\tilde{v}^\beta+\tilde{g}^{\alpha\beta}P\right]
\end{equation*}

\begin{equation*}
= r\left[\frac{\partial}{\partial \xi^\alpha}(\rho \sqrt{g}v^\alpha V^3) + 2\rho\sqrt{g}(V^3)^2\right] + 2r\sqrt{g}P - rg_{\alpha\beta}\rho\sqrt{g}v^\alpha v^\beta - rg_{\alpha\beta}\sqrt{g}g^{\alpha\beta}P
\end{equation*}

the terms involving $P$ cancel leaving:

\begin{equation*}
= r\left[\frac{\partial}{\partial \xi^\alpha}(\rho \sqrt{g}v^\alpha V^3) + 2\rho\sqrt{g}(V^3)^2 - \rho\sqrt{g}q^2_c \right]
\end{equation*}

Dividing by $r$ gives the conservation form for the radial momentum component:

\begin{equation*}
\frac{\partial}{\partial \xi^\alpha}(\rho \sqrt{g}v^\alpha V^3) + 2\rho\sqrt{g}(V^3)^2 - \rho\sqrt{g}q^2_c = 0
\end{equation*}

which can be rearranged as:

\begin{equation*}
V^3\left[\frac{\partial}{\partial \xi^\alpha}(\rho \sqrt{g}v^\alpha ) + 2\rho\sqrt{g}V^3\right] + \rho \sqrt{g}v^\alpha\frac{\partial V^3}{\partial \xi^\alpha} - \rho\sqrt{g}q_c^2 = 0
\end{equation*}

The term in brackets is the mass continuity equation and is thus equal to zero. Thus:

\begin{equation}
\rho \sqrt{g}v^\alpha\frac{\partial V^3}{\partial \xi^\alpha} - \rho\sqrt{g}q^2_c = 0
\end{equation}

and cancelling $\rho\sqrt{g}$:

\begin{equation}
v^\alpha\frac{\partial V^3}{\partial \xi^\alpha} - q^2_c = 0
\end{equation}

This closes the system.

\section{Comparison to spherical coordinates }\label{sphcoords}

In order to validate this derivation, we verify that the equations reduce to a form given in past literature when the geometric relationships of spherical coordinates are plugged in. The conical continuity and momentum equations were presented in \cite{NASA_con} using traditional spherical coordinates. With $\phi$ as the zenith angle and $\theta$ as the azimuthal angle, these equations are:

\begin{subequations}\label{spherical}
\begin{gather}
\frac{\partial}{\partial \phi}\left( \rho v^\phi\sin\phi \right) + \frac{\partial }{\partial \theta}\left(\rho v^\theta\sin\phi \right) + 2\rho V^r\sin\phi = 0 \label{sph_mass}
\\
v^\phi\frac{\partial v^\phi}{\partial \phi} + v^\theta\frac{\partial v^\phi}{\partial \theta} + \frac{1}{\rho}\frac{\partial P}{\partial \phi} + v^\phi V^r - (v^\theta)^2\sin\phi\cos\phi = 0 \label{sph_mom_phi}
\\
v^\phi\frac{\partial }{\partial \phi}(v^\theta\sin\phi) + v^\theta \sin\phi\frac{\partial v^\theta}{\partial \theta} + \frac{1}{\rho\sin\phi}\frac{\partial P}{\partial \theta} + v^\theta V^r\sin\phi + v^\theta v^\phi\cos\phi = 0 \label{sph_mom_th}
\\
v^\phi\frac{\partial V^r}{\partial \phi} + v^\theta\frac{\partial V^r}{\partial \theta} - (v^\phi)^2 - (v^\theta)^2\sin^2\phi = 0 \label{sph_mom_r}
\end{gather}
\end{subequations}

By letting $\phi\sim1$ and $\theta\sim2$, the metric tensor for spherical coordinates on the surface is given by:

\begin{equation}
    g_{\alpha\beta} = \left[ \begin{matrix} 1 & 0 \\ 0 & \sin^2\phi \end{matrix} \right]
\end{equation}

therefore:

\begin{equation}
    g^{\alpha\beta} = \left[ \begin{matrix} 1 & 0 \\ 0 & \frac{1}{\sin^2\phi} \end{matrix} \right]
\end{equation}

and

\begin{equation}
    g=\sin^2\phi
\end{equation}

The Christoffel symbols are given by:

\begin{equation}
    \overset{(g)}{\Gamma}\indices{_\gamma^\phi_\nu} = \left[ \begin{matrix} 0 & 0 \\ 0 & -\sin\phi\cos\phi \end{matrix} \right]
\end{equation}

\begin{equation}
    \overset{(g)}{\Gamma}\indices{_\gamma^\theta_\nu} = \left[ \begin{matrix} 0 & \frac{\cos\phi}{\sin\phi} \\ \frac{\cos\phi}{\sin\phi} & 0 \end{matrix} \right]
\end{equation}

With these, Equation \eqref{sph_mass} can be obtained from Equation \eqref{mass} by direct substitution. To get Equations \eqref{sph_mom_phi} and \eqref{sph_mom_th}, one must first subtract Equation \eqref{mass} from Equation \eqref{mom}, and then divide the resulting $v^1$ equation by $\rho\sqrt{g}$ and the resulting $v^2$ equation by just $\rho$. The rest follows from direct substitution of the expressions for the metric and its determinant and the Christoffel symbols and some manipulation of the derivatives.  Lastly, Equation \eqref{sph_mom_r} comes through direct substitution of the expression for $q_c$ into equation \eqref{mom_r}.

System \ref{TheEq} is thus consistent with traditional spherical coordinates which demonstrates that the new system of equations was derived correctly. We emphasize here again that the new system presented here is superior to a formulation in traditional spherical coordinates in that it is easily adapted to general coordinate systems where there may not be such simple or established geometric relationships.

\section{Elliptic-Hyperbolic property}\label{hyperbolic}

A general 1st-order system of $m$ differential equations in $n$ spatial dimensions has the following form:

\begin{equation}\label{general}
U_t + \sum_{i=1}^{n}\bar{A}^iU_{x^i}+\bar{S}=0
\end{equation}

Where $U:\mathbb{R}^n\rightarrow\mathbb{R}^m$ is a column vector of the dependent variables and each $\bar{A}^i$ is an $m$ by $m$ matrix that can in general depend on $U$ and $\pmb{x}$. S is a column vector of source terms. 

\begin{definition}
A system of the form \eqref{general} is said to be strictly hyperbolic if $\forall \pmb{w}\in \mathbb{R}^n, |\pmb{w}|=1$, the eigenvalues of $\bar{A}_w=\sum_{i=1}^{n}w_i\bar{A}^i$ are real and distinct. If they are all real, but not all distinct, the system is non-strictly hyperbolic. If any of the eigenvalues are complex then the system is said to be elliptic \cite{evansPDE,lax2006hyperbolic,renardy,CHEN20051}.
\end{definition}

%attribute this definition to Evans Chen and Lax

\section{Eigenvalues}\label{valsvects}

In equation \eqref{TheEq}, the dependent variables are:

\begin{equation}
U = \left[ \begin{matrix} \rho \\ v^1 \\ v^2 \\ V^3 \\ e \end{matrix} \right]
\end{equation}

After using the product rule and/or chain rule to expand all the derivatives until they are in terms of derivatives of individual dependent variables, the system has the form:

\begin{equation}\label{general_steady}
\sum_{\alpha=1}^{2}A^\alpha U_{\xi^\alpha}+S=0
\end{equation}

where $A^\alpha$ is the jacobian matrix of:

\begin{equation*}
\left[ \begin{matrix} \rho\sqrt{g}v^\alpha \\ \sqrt{g}[\rho v^1v^\alpha + g^{1\alpha}P] \\ 
    \sqrt{g}[\rho v^2v^\alpha + g^{2\alpha}P] \\ \rho\sqrt{g}V^3v^\alpha \\ \sqrt{g}[\rho E+P]v^\alpha \end{matrix} \right]
\end{equation*}

given by:

\begin{equation}
    A^1 = \sqrt{g}\left[
\begin{matrix}
   v^1 & \rho & 0 & 0 & 0 \\
   (v^1)^2+ g^{11}P_\rho & 2 \rho v^1 &  0 &  0 &  g^{11} P_e \\
 v^1 v^2 +  g^{12} P_\rho   & \rho v^2  & \rho v^1  &  0 &  g^{12} P_e \\
   v^1 V^3 & \rho V^3 & 0 & \rho v^1 & 0 \\
   v^1 \left(E+P_\rho\right) & \rho (g_{1\beta}v^\beta) v^1 + (\rho E+P) & \rho (g_{2\beta}v^\beta) v^1 & \rho V^3v^1 & v^1 \left(\rho+P_e\right) \\
\end{matrix}
\right]
\end{equation}

and

\begin{equation}
    A^2 = \sqrt{g}\left[
\begin{matrix}
   v^2 & 0 & \rho & 0 & 0 \\
   v^1 v^2+ g^{21}P_\rho & \rho v^2 & \rho v^1  &  0 &  g^{21} P_e \\
 (v^2)^2 + g^{22} P_\rho & 0 & 2\rho v^2 &  0 &  g^{22} P_e \\
   v^2 V^3 & 0 & \rho V^3 & \rho v^2 & 0 \\
   v^2 \left(E+P_\rho\right) & \rho(g_{1\beta}v^\beta)v^2 & \rho(g_{2\beta}v^\beta)v^2 + (\rho E+P) & \rho V^3v^2 & v^2 \left(\rho+P_e\right) \\
\end{matrix}
\right]
\end{equation}

For hyperbolicity to be assessed, a spatial variable must be chosen to be treated as time-like. Without loss of generality, $\xi^1$ is chosen. System \eqref{general_steady} is then multiplied by the inverse of $A^1$ giving:

\begin{equation}
U_{\xi^1}+\bar{A}U_{\xi^2}+\bar{S}=0
\end{equation}

where $\bar{A} = (A^1)^{-1}A^2$ and $\bar{S}=(A^1)^{-1}S$. There is only one matrix left, so we simply take $w=1$ and $\bar{A}_w = \bar{A}$. The eigenvalues of $\bar{A}$ were computed using Wolfram Mathematica \cite{Mathematica} and are given in Equation \eqref{FullEvals}.

\begin{equation}\label{FullEvals}
\lambda(\bar{A})=\frac{v^2}{v^1},\frac{v^2}{v^1},\frac{v^2}{v^1},\\
\frac{v^1v^2 - c^2g^{12}\mp \frac{c}{\sqrt{g}}\sqrt{q_c^2-c^2}}{(v^1)^2-g^{11}c^2}
\end{equation}

where the speed of sound, $c$, is given by:

\begin{equation}
    c = \frac{\sqrt{PP_e + \rho^2 P_\rho} }{\rho}
\end{equation}

\begin{remark}
We note that in the event that $P$ is a function of $\rho$ only this expression reduces to the expression:

\begin{equation}
    c = \sqrt{\frac{\partial P}{\partial \rho}}
\end{equation}

and for an ideal gas with $P=(\gamma-1)\rho e$ this gives:

\begin{equation}
    c = \sqrt{\frac{\gamma P}{\rho}}
\end{equation}
\end{remark}

The first three eigenvalues coincide so the system cannot be strictly hyperbolic. The last two eigenvalues will become complex if the magnitude of the crossflow velocity is less than the speed of sound. This means that the type is:

\begin{equation}
\text{type}=\begin{cases}
	\text{hyperbolic} & q_c > c \\ 
    \text{elliptic} & q_c < c \\ 
\end{cases}
\end{equation}

Which is a result analogous to steady Euler in the Cartesian setting.

\subsection{Pseudo-Time Dependency}

The time derivative terms in Equation \eqref{transient} are not compatible with the conical assumption because the $r$ dependency fails to disappear. However, it is often convenient to solve a steady problem numerically by marching in time until the solution stabilizes. For that purpose one could reinsert the time derivatives with the appropriate metrics and treat the problem as unsteady. 

For this nonphysical problem, the form is:

\begin{equation}
A_0U_t + \sum_{\alpha=1}^{2}A^\alpha U_{\xi^\alpha}=0
\end{equation}

And after multiplying by the inverse of $A_0$:

\begin{equation}
U_t + \sum_{\alpha=1}^{2}\bar{A}^\alpha U_{\xi^\alpha}=0
\end{equation}

defining $\bar{A}^\alpha = A_0^{-1}A^\alpha$. Since the system is given for the contravariant components of the velocity we take $\pmb{w}$ to be a covariant vector such that $g^{\alpha\beta}w_\alpha w_\beta = 1$, and form the linear combination $\bar{A}_w=\sum_{i=1}^{n}w_i\bar{A}^i$. The eigenvalues are:

\begin{equation}
    \lambda(\bar{A}_w) = \pmb{v}\cdot\pmb{w}, \pmb{v}\cdot\pmb{w}, \pmb{v}\cdot\pmb{w}, \pmb{v}\cdot\pmb{w}\pm c 
\end{equation}

which are the same as for general unsteady Euler. Therefore, this system is everywhere non-strictly hyperbolic which differs from the steady case.

\section{Comparison to potential case}\label{potential}

The case of potential flow was examined by Sritharan in \cite{sriThesis} which is a similar, but simpler problem compared to that considered here, subject to the additional assumption that the flow is irrotational. It is thus relevant to see how the wave speeds compare. For that case, the mass continuity equation was identical to equation \eqref{mass}, the momentum equation was replaced with the equation $\pmb{V}=\nabla\psi$ where $\psi$ is some scalar function, and an isentropic energy equation was used to relate density to velocity and the free stream conditions. For the conical assumption to hold, the function whose gradient is the velocity has the form $\psi=rF(\xi^1,\xi^2)$. Therefore:

\begin{equation}
    V_i = \frac{\partial }{\partial x^i}(rF)
\end{equation}

which in turn gives for the scaled surface components:

\begin{equation}
    v_\alpha = \frac{\partial F}{\partial \xi^\alpha}
\end{equation}

Note that this immediately satisfies Equation \eqref{mom_r}. With $V_3=V^3=F$, the first term on the left becomes:

\begin{equation}
v^\alpha\frac{\partial F}{\partial \xi^\alpha} = v^\alpha v_\alpha = q^2_c
\end{equation}

which cancels with the second term on the left thus satisfying the equation.

With this expression for the velocity, equation \eqref{mass} becomes:

\begin{equation}
    \frac{\partial}{\partial\xi^\beta}\left( \rho\sqrt{g}g^{\beta\alpha}\frac{\partial F}{\partial\xi^\alpha} \right) + 2\rho\sqrt{g}F = 0
\end{equation}

and after plugging in the energy equation to relate density to the velocity, the result was the second order equation:

\begin{equation}
    \rho\sqrt{g}\left[\left(g^{\alpha\beta}-\frac{v^\alpha v^\beta}{c^2}\right)\frac{\partial^2F}{\partial\xi^\alpha\partial\xi^\beta} + H\left(F,\frac{\partial F}{\partial\xi^\alpha},g^{\alpha\beta}\right)\right]=0
\end{equation}

where $H$ is the collection of terms involving lower order derivatives. The second order portion can be converted to the following first order system:

\begin{equation}
    \begin{bmatrix} g^{11}-\frac{(v^1)^2}{c^2} & g^{12}-\frac{v^1v^2}{c^2} \\ 0 & 1  \end{bmatrix} \begin{bmatrix} v_1 \\ v_2  \end{bmatrix}_{\xi^1} + \begin{bmatrix} g^{12}-\frac{v^1v^2}{c^2} & g^{22}-\frac{(v^2)^2}{c^2} \\ -1 & 0  \end{bmatrix} \begin{bmatrix} v_1 \\ v_2  \end{bmatrix}_{\xi^2}
\end{equation}

After multiplying by the inverse of the leading matrix, the eigenvalues of the matrix on the second term come out to be:

\begin{equation}
\lambda=\frac{v^1v^2 - c^2g^{12}\pm \frac{c}{\sqrt{g}}\sqrt{q_c^2-c^2}}{(v^1)^2-g^{11}c^2}
\end{equation}

which are the last two eigenvalues of the full system. The potential case thus has the same type and some of the same wave speeds as the general flow case, demonstrating the consistency between the two cases and further validating the derivation presented in this paper.

\section{Conclusion}

We have thus systematically derived a flexible system of equations which describe inviscid flow past a cone of arbitrary cross section. Under the assumption of conical invariance, the compressible Euler equations reduce to a system which is defined on the unit sphere. The mathematical formulation of the system given in this article is superior to systems provided in the past as it is stated for a general curved coordinate systems and is ready to interface with modern grid generation and numerical solution methods. It has also been demonstrated that the type can be either hyperbolic or elliptic depending on whether the crossflow Mach number on the sphere is supersonic or subsonic.

\section*{Acknowledgement}

This research was supported in part by an appointment to the Student Research Participation Program at the U.S. Air Force Institute of Technology administered by the Oak Ridge Institute for Science and Education through an interagency agreement between the U.S. Department of Energy and USAFIT.

%%%%%%%%%% Insert bibliography here %%%%%%%%%%%%%%

\bibliographystyle{amsplain}
\bibliography{references.bib}

\end{document}